\documentclass[12pt,a4paper,twoside]{article}
\usepackage{amssymb,amsmath,amscd,amsthm}
\usepackage{epic}
\usepackage{eepic}

\input{epsf}

\newcommand{\pro}[2]{\langle #1, #2 \rangle}

\def\P{{\mathbb P}}
\def\Q{{\mathbb Q}}
\def\R{{\mathbb R}}
\def\Z{{\mathbb Z}}

\def\F{{\cal F}}

\def\J{{\cal R}}

\def\V{{{\cal V}}}

\def\Hom{{\textup{Hom}}}

\def\NR{{N_{\R}}}

\def\MR{{M_{\R}}}

\def\rand{\partial}
\def\randp{{\rand P}}

\newtheorem*{theorem*}{Theorem}
\newtheorem*{corollary*}{Corollary}
\newtheorem{theorem}{Theorem}[section]
\newtheorem{corollary}[theorem]{Corollary}

\newtheorem{lemma}[theorem]{Lemma}

\theoremstyle{definition}
\newtheorem{definition}[theorem]{Definition}

\title{Lattices generated by skeletons \\of reflexive polytopes}

\begin{document}

\date{}

\author{Christian Haase\\{\small Mathematics Department, Duke University, }\\
{\small Durham, NC 27708, USA}\\
{\small email: haase@math.duke.edu}\\\vspace{0.01cm}\\
Benjamin Nill\\{\small Mathematisches Institut, Universit\"at T\"ubingen}\\
{\small Auf der Morgenstelle 10, 72076 T\"ubingen, Germany}\\
{\small email: nill@algebra.mathematik.uni-tuebingen.de}}

\maketitle

\abstract{Lattices generated by lattice points in skeletons of
  reflexive polytopes are essential in determining the fundamental
  group and integral cohomology of Calabi-Yau hypersurfaces.
  Here we prove that the lattice generated by all lattice points in a
  reflexive polytope is already generated by lattice points in
  codimension two faces. This answers a question of J. Morgan.
}

\section{Introduction and main result}
Since its introduction by Batyrev in \cite{Bat94} dual pairs of
reflexive polytopes have been used to successfully construct mirror
pairs of smooth Calabi-Yau varieties as resolutions of non-degenerate
anticanonical hypersurfaces in Gorenstein toric Fano varieties,
see, e.g., \cite{CK99}. 

Recall that a reflexive polytope is an $n$-dimensional lattice
polytope $P \subseteq \MR = M \otimes_\Z \R$ for a lattice $M \cong
\Z^n$ such that $P$ contains the origin in its interior and the dual
polytope $P^*$
is also a lattice polytope with respect to the dual lattice $N = M^* =
\Hom_\Z(M,\Z)$. There is an associated pair $\P_{\Sigma_P,M}$,
$\P_{\Sigma_{P^*},N}$ of Gorenstein toric Fano varieties, where
$\Sigma_P$ is just the fan of cones over the faces of $P$, see, e.g.,
\cite{Bat94,Nil05}.


For a given reflexive polytope $P \subseteq \MR$ there exist finitely
many choices of sublattices of $M$ such that $P$ is reflexive with
respect to this lattice. Obviously the minimal one of these is the
lattice $\Lambda_0$ generated by the vertices of $P$. More generally
we define \[\Lambda_k := \text{lattice generated by all lattice
  points in the $k$-skeleton of $P$,}\]
where the $k$-skeleton is the union of $k$-dimensional faces. 

Due to \cite[Lemma 1.17]{Nil05} and \cite[3.2]{Ful93} we see
that for $k = 0,1,2$ the quotient group $M / \Lambda_k$ equals the
fundamental group $\pi_1$ of the union of all $\leq k+1$-codimensional
torus orbits in $\P_{\Sigma_P,M}$. 

Returning to the relevance of these lattices in mirror symmetry
Batyrev and Kreuzer \cite[Cor.1.9]{BK05} recently proved that for 
$n = 3,4$ the quotient group $M / \Lambda_{n-2}$ is precisely the
fundamental group $\pi_1(X)$ of a projective crepant resolution $X$ of
an $(n-1)$-dimensional non-degenerate anticanonical Calabi-Yau
hypersurface in $\P_{\Sigma_P,M}$. Moreover for $n=4$ they showed that
the torsion group of $H^2(X;\Z)$ equals $\Hom(M / \Lambda_{n-2},
\Q/\Z)$ and the torsion group of $H^3(X;\Z)$ is isomorphic to $\Hom(\wedge^2 
M / M \wedge \Lambda_{n-3}, \Q/\Z)$, cf. \cite[Cor.3.9]{BK05}.
Mirror symmetry should exchange these torsion parts of integral
cohomology of $X$ and a mirror $X^*$. Hence, this yields the
surprising isomorphism $M / \Lambda_2 \cong \wedge^2 N / N \wedge
\Lambda^*_1$, where $\Lambda^*_1$ 
is the sublattice of $N$ generated
by lattice points in edges of $P^*$. Using the classification of Kreuzer and
Skarke \cite{KS00} they confirmed this conjecture.


These results and open questions motivate further investigation of the
lattices $\Lambda_k$ for an $n$-dimensional reflexive polytope
$P$. Here we first note the following two observations:
\begin{itemize}
\item Since there are no non-zero lattice points in the interior of
  $P$, we have $\Lambda_{n-1} = \Lambda_n$. 
\item If there exists a crepant toric resolution of $\P_{\Sigma_P,M}$,
  then the boundary $\randp$ contains a lattice basis, hence
  $\Lambda_{n-1} = M$. This holds for $n=2,3$.
\end{itemize}

For $n=2$ there are precisely three isomorphism classes of reflexive
polytopes with $\Lambda_{n-2} = \Lambda_0 \not= M$, here $\Lambda_0$ has index $2$, $2$, and $3$:
\begin{figure}[htbp]
  \centering
  \setlength{\unitlength}{0.00016667in}
\begingroup\makeatletter\ifx\SetFigFont\undefined%
\gdef\SetFigFont#1#2#3#4#5{%
  \reset@font\fontsize{#1}{#2pt}%
  \fontfamily{#3}\fontseries{#4}\fontshape{#5}%
  \selectfont}%
\fi\endgroup%
{\renewcommand{\dashlinestretch}{30}
\begin{picture}(13366,3781)(0,-10)
\path(83,83)(2483,2483)(4883,83)(83,83)
\path(6083,2483)(6083,83)(8483,83)
	(8483,2483)(6083,2483)
\path(9683,3683)(9683,83)(13283,83)(9683,3683)
\put(13283,3683){\blacken\ellipse{150}{150}}
\put(13283,3683){\ellipse{150}{150}}
\put(12083,3683){\blacken\ellipse{150}{150}}
\put(12083,3683){\ellipse{150}{150}}
\put(10883,3683){\blacken\ellipse{150}{150}}
\put(10883,3683){\ellipse{150}{150}}
\put(9683,3683){\blacken\ellipse{150}{150}}
\put(9683,3683){\ellipse{150}{150}}
\put(13283,2483){\blacken\ellipse{150}{150}}
\put(13283,2483){\ellipse{150}{150}}
\put(12083,2483){\blacken\ellipse{150}{150}}
\put(12083,2483){\ellipse{150}{150}}
\put(10883,2483){\blacken\ellipse{150}{150}}
\put(10883,2483){\ellipse{150}{150}}
\put(9683,2483){\blacken\ellipse{150}{150}}
\put(9683,2483){\ellipse{150}{150}}
\put(13283,1283){\blacken\ellipse{150}{150}}
\put(13283,1283){\ellipse{150}{150}}
\put(12083,1283){\blacken\ellipse{150}{150}}
\put(12083,1283){\ellipse{150}{150}}
\put(10883,1283){\blacken\ellipse{150}{150}}
\put(10883,1283){\ellipse{150}{150}}
\put(9683,1283){\blacken\ellipse{150}{150}}
\put(9683,1283){\ellipse{150}{150}}
\put(13283,83){\blacken\ellipse{150}{150}}
\put(13283,83){\ellipse{150}{150}}
\put(12083,83){\blacken\ellipse{150}{150}}
\put(12083,83){\ellipse{150}{150}}
\put(10883,83){\blacken\ellipse{150}{150}}
\put(10883,83){\ellipse{150}{150}}
\put(9683,83){\blacken\ellipse{150}{150}}
\put(9683,83){\ellipse{150}{150}}
\put(8483,83){\blacken\ellipse{150}{150}}
\put(8483,83){\ellipse{150}{150}}
\put(7283,83){\blacken\ellipse{150}{150}}
\put(7283,83){\ellipse{150}{150}}
\put(6083,83){\blacken\ellipse{150}{150}}
\put(6083,83){\ellipse{150}{150}}
\put(4883,83){\blacken\ellipse{150}{150}}
\put(4883,83){\ellipse{150}{150}}
\put(3683,83){\blacken\ellipse{150}{150}}
\put(3683,83){\ellipse{150}{150}}
\put(2483,83){\blacken\ellipse{150}{150}}
\put(2483,83){\ellipse{150}{150}}
\put(1283,83){\blacken\ellipse{150}{150}}
\put(1283,83){\ellipse{150}{150}}
\put(83,83){\blacken\ellipse{150}{150}}
\put(83,83){\ellipse{150}{150}}
\put(8483,1283){\blacken\ellipse{150}{150}}
\put(8483,1283){\ellipse{150}{150}}
\put(7283,1283){\blacken\ellipse{150}{150}}
\put(7283,1283){\ellipse{150}{150}}
\put(6083,1283){\blacken\ellipse{150}{150}}
\put(6083,1283){\ellipse{150}{150}}
\put(4883,1283){\blacken\ellipse{150}{150}}
\put(4883,1283){\ellipse{150}{150}}
\put(3683,1283){\blacken\ellipse{150}{150}}
\put(3683,1283){\ellipse{150}{150}}
\put(2483,1283){\blacken\ellipse{150}{150}}
\put(2483,1283){\ellipse{150}{150}}
\put(1283,1283){\blacken\ellipse{150}{150}}
\put(1283,1283){\ellipse{150}{150}}
\put(83,1283){\blacken\ellipse{150}{150}}
\put(83,1283){\ellipse{150}{150}}
\put(8483,2483){\blacken\ellipse{150}{150}}
\put(8483,2483){\ellipse{150}{150}}
\put(7283,2483){\blacken\ellipse{150}{150}}
\put(7283,2483){\ellipse{150}{150}}
\put(6083,2483){\blacken\ellipse{150}{150}}
\put(6083,2483){\ellipse{150}{150}}
\put(4883,2483){\blacken\ellipse{150}{150}}
\put(4883,2483){\ellipse{150}{150}}
\put(3683,2483){\blacken\ellipse{150}{150}}
\put(3683,2483){\ellipse{150}{150}}
\put(2483,2483){\blacken\ellipse{150}{150}}
\put(2483,2483){\ellipse{150}{150}}
\put(1283,2483){\blacken\ellipse{150}{150}}
\put(1283,2483){\ellipse{150}{150}}
\put(83,2483){\blacken\ellipse{150}{150}}
\put(83,2483){\ellipse{150}{150}}
\put(8483,3683){\blacken\ellipse{150}{150}}
\put(8483,3683){\ellipse{150}{150}}
\put(7283,3683){\blacken\ellipse{150}{150}}
\put(7283,3683){\ellipse{150}{150}}
\put(6083,3683){\blacken\ellipse{150}{150}}
\put(6083,3683){\ellipse{150}{150}}
\put(4883,3683){\blacken\ellipse{150}{150}}
\put(4883,3683){\ellipse{150}{150}}
\put(3683,3683){\blacken\ellipse{150}{150}}
\put(3683,3683){\ellipse{150}{150}}
\put(2483,3683){\blacken\ellipse{150}{150}}
\put(2483,3683){\ellipse{150}{150}}
\put(1283,3683){\blacken\ellipse{150}{150}}
\put(1283,3683){\ellipse{150}{150}}
\put(83,3683){\blacken\ellipse{150}{150}}
\put(83,3683){\ellipse{150}{150}}
\end{picture}
}
  \caption{Reflexive polygons with $\Lambda_{n-2} \not= M$}
  \label{fig:counterEG}
\end{figure}
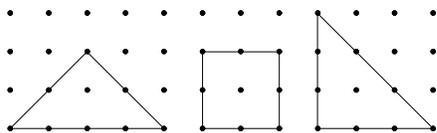

As noted in \cite[Remark 1.10]{BK05} for $n=3$ a non-degenerate
Calabi-Yau hypersurface $X$ arising from a crepant resolution of
$\P_{\Sigma_P,M}$ is a smooth $K3$-surface, hence simply-connected. So
the fundamental group $\pi_1(X) \cong M / \Lambda_{n-2}$ is trivial,
i.e., $\Lambda_{n-2} = \Lambda_{n-1}$. The authors did not know of a
combinatorial proof of this result in the literature.


At PCMI 2004, John Morgan (for a mirror-symmetric conjecture concerning integral variations of Hodge structure 
see Conjecture 2.9 in \cite{DM05}) 
asked the first author for a combinatorial proof that
$\Lambda_{n-2} = \Lambda_{n-1}$ would also hold for $n=4$. This was
confirmed by the classification of all four-dimensional reflexive
polytopes \cite{KS00}. John Morgan, and independently Batyrev and
Kreuzer \cite{BK05}, have listed all $16$ isomorphism classes of
four-dimensional reflexive polytopes with $\Lambda_{n-2} \not=
M$. However in contrast to $n=3$ not even an algebro-geometric proof
seemed to be known.


The goal of this article is to provide a purely convex-geometric proof
valid in arbitrary dimension $n > 2$. (Figure \ref{fig:counterEG}
shows all ``counter-examples'' for $n=2$).

\begin{theorem} \label{thm:main}
  If $n \geq 3$, then $\Lambda_{n-2} = \Lambda_{n-1}$.
\end{theorem}

The main idea of the proof is to show that any lattice point inside a
facet can be obtained from lattice points on lower dimensional
faces. For this we use a partial addition on the set of lattice points
in $P$, observed by the second author in \cite[Prop.4.1]{Nil05}. In
\cite{Nil05} that method has been applied to investigate the set of
lattice points inside the facets, since this is precisely the set of
Demazure roots for the automorphism group of the corresponding
Gorenstein toric Fano variety $\P_{\Sigma_{P^*},N}$.

\bigskip

\textbf{Acknowledgments:}
Work of the first author was partially supported by NSF grant
DMS-0200740. The second author would like to thank Duke University 
for hospitality and support.

\section{Proof of theorem}

Let $P \subseteq \MR$ be an $n$-dimensional reflexive polytope. We denote by $\pro{\cdot}{\cdot}$ the 
non-degenerate symmetric pairing of the dual lattices $M,N$.

For the proof we will use the following notions.
\begin{definition}
\begin{itemize}
\item[]
\item We denote by $\J$ the set of {\em Demazure roots} of $P$, i.e.,
  the set of lattice points in the interior of facets.

\item For $x \in \J$ we denote by $\F_x$ the unique facet of $P$
  containing $x$, and by $\eta_x$ the unique inner normal of $\F_x$,
  i.e., $\eta_x = \eta_{\F_x} \in N$ s.t. $\pro{\eta_x}{\F_x} = -1$.

\item A pair of roots $x,y \in \J$ is called {\em orthogonal}, if
  $\pro{\eta_x}{y} = 0 = \pro{\eta_y}{x}$. In particular $\eta_{-x}
  \not= \eta_y \not= \eta_x \not= \eta_{-y}$. In this case $x + y \in
  M \cap \F_x \cap \F_y$ (see the next result).
\end{itemize}
\end{definition}
The following lemma is the engine running our proof.
\begin{lemma}[{\cite[Lemma 3.5 \& Corollary 3.6]{Nil04}}] \label{lemma:nill}
  For $x,y \in \randp \cap M$ exactly one of the following is true.
  Either $x$ and $y$ belong to a common facet, or $x+y=0$, or $x+y \in
  \randp$. In the last case, there is a unique $z(x,y) = ax + by \in
  \randp$ with $a,b \in \Z_{>0}$. We have $a=1$ or $b=1$, and 
$x$ and $z(x,y)$, respectively $y$ and $z(x,y)$, belong to a common facet.
  Furthermore, if $x \in \J$, 
  then
  $a = \pro{\eta_x}{y} + 1$.
\end{lemma}
The custom tailored version for our purposes reads as
follows.\footnote{Of course, this corollary talks about the empty set
  as we aim to show.} We abbreviate $\Lambda := \Lambda_{n-2}$.
\begin{corollary} \label{tool}
  Let $x \in \J \setminus \Lambda$, $y \in \randp \cap \Lambda$ with
  $y \not\in \F_x$. Then $\pro{\eta_x}{y} \geq 1$, and $p(x,y) :=
  \pro{\eta_x}{y} x + y \in \J \setminus \Lambda$ is a root orthogonal
  to $x$.
\end{corollary}
\begin{figure}[htbp]
  \centering
  \setlength{\unitlength}{0.00016667in}
\begingroup\makeatletter\ifx\SetFigFont\undefined%
\gdef\SetFigFont#1#2#3#4#5{%
  \reset@font\fontsize{#1}{#2pt}%
  \fontfamily{#3}\fontseries{#4}\fontshape{#5}%
  \selectfont}%
\fi\endgroup%
{\renewcommand{\dashlinestretch}{30}
\begin{picture}(7224,5540)(0,-10)
\put(5187,113){\blacken\ellipse{212}{212}}
\put(5187,113){\ellipse{212}{212}}
\put(5187,2213){\blacken\ellipse{212}{212}}
\put(5187,2213){\ellipse{212}{212}}
\put(5187,4313){\blacken\ellipse{212}{212}}
\put(5187,4313){\ellipse{212}{212}}
\put(3012,4313){\blacken\ellipse{212}{212}}
\put(3012,4313){\ellipse{212}{212}}
\path(12,3713)(4212,3713)(4212,113)
\path(3012,5513)(7212,5513)(7212,2663)
\path(4212,3713)(7212,5513)
\path(7212,1913)(7212,1613)
\dottedline{150}(3087,4313)(5187,4313)(5187,113)
\put(5487,188){\makebox(0,0)[lb]{\smash{{{\SetFigFont{7}{8.4}{\rmdefault}{\mddefault}{\updefault}$y$}}}}}
\put(2937,4538){\makebox(0,0)[lb]{\smash{{{\SetFigFont{7}{8.4}{\rmdefault}{\mddefault}{\updefault}$x$}}}}}
\put(5487,2063){\makebox(0,0)[lb]{\smash{{{\SetFigFont{7}{8.4}{\rmdefault}{\mddefault}{\updefault}$p(x,y)=\pro{\eta_x}{y}x+y$}}}}}
\put(4212,4763){\makebox(0,0)[lb]{\smash{{{\SetFigFont{7}{8.4}{\rmdefault}{\mddefault}{\updefault}$\F_x$}}}}}
\end{picture}
}
  \caption{$p(x,y)$}
  \label{fig:pxy}
\end{figure}
\begin{proof}[Proof of Theorem~\ref{thm:main}]\
Let $n \geq 3$, and let $x_1 \in \J$ arbitrary. We have to show $x_1
\in \Lambda$. Assume not, that is, $x_1 \not\in \Lambda$.

Let $y_0 \in \V(P) \setminus \F_{x_1}$ be a vertex of $P$ outside
$\F_{x_1}$. Corollary~\ref{tool} provides $x_2 := p(x_1,y_0) \in \J
\setminus \Lambda$.
The points $x_1$ and $x_2$ span a $2$-dimensional linear space $L$.
Either all vertices of $P$ belong to $(\F_{x_1} \cap \F_{x_2}) \cup L$
or not.
\begin{itemize}
\item
\underline{$\V(P) \subset (\F_{x_1} \cap \F_{x_2}) \cup L$ :} The
intersection $Q = P \cap L$ is a reflexive lattice polygon since any
vertex of $Q$ that is not on $\F_{x_1} \cap \F_{x_2}$ is a vertex of
$P$, and the other vertex $x_1+x_2$ is integral, too. The part
of $\partial Q$ outside $\F_{x_1} \cup \F_{x_2}$ belongs to the
$1$-skeleton of $P$. The endpoints of a primitive segment form a
lattice basis for $L \cap M \ni x_1$. So $x_1 \in \Lambda_1 \subseteq \Lambda$, a contradiction. 
(This is where we use $n \ge 3$ !)
\begin{figure}[htbp]
  \centering
  \setlength{\unitlength}{0.00016667in}
\begingroup\makeatletter\ifx\SetFigFont\undefined%
\gdef\SetFigFont#1#2#3#4#5{%
  \reset@font\fontsize{#1}{#2pt}%
  \fontfamily{#3}\fontseries{#4}\fontshape{#5}%
  \selectfont}%
\fi\endgroup%
{\renewcommand{\dashlinestretch}{30}
\begin{picture}(6849,5615)(0,-10)
\put(4212,113){\blacken\ellipse{212}{212}}
\put(4212,113){\ellipse{212}{212}}
\put(4212,2213){\blacken\ellipse{212}{212}}
\put(4212,2213){\ellipse{212}{212}}
\put(2112,4313){\blacken\ellipse{212}{212}}
\put(2112,4313){\ellipse{212}{212}}
\path(12,4313)(3237,3713)(4212,113)
\path(12,4313)(6237,5513)
\path(5712,4088)(6237,5513)
\path(5157,2618)(5622,3863)
\path(6837,4388)(6835,4386)(6832,4382)
	(6825,4375)(6815,4364)(6801,4349)
	(6783,4331)(6761,4309)(6737,4284)
	(6711,4258)(6683,4231)(6655,4203)
	(6626,4176)(6597,4151)(6570,4126)
	(6543,4103)(6516,4082)(6491,4063)
	(6467,4045)(6443,4029)(6419,4015)
	(6396,4002)(6372,3991)(6348,3981)
	(6324,3971)(6300,3963)(6276,3956)
	(6252,3950)(6226,3944)(6200,3939)
	(6174,3934)(6146,3931)(6117,3928)
	(6088,3925)(6058,3924)(6027,3923)
	(5996,3923)(5964,3924)(5932,3926)
	(5899,3928)(5867,3932)(5835,3936)
	(5803,3940)(5772,3946)(5741,3952)
	(5711,3959)(5682,3967)(5653,3975)
	(5625,3984)(5599,3994)(5573,4004)
	(5547,4015)(5523,4026)(5500,4038)
	(5476,4051)(5454,4064)(5432,4078)
	(5410,4094)(5388,4111)(5366,4129)
	(5344,4148)(5322,4170)(5299,4193)
	(5275,4218)(5251,4245)(5225,4274)
	(5200,4304)(5173,4336)(5147,4369)
	(5120,4403)(5094,4436)(5069,4468)
	(5046,4499)(5025,4527)(5007,4551)
	(4992,4572)(4980,4588)(4962,4613)
\path(5150.925,4453.290)(4962.000,4613.000)(5053.541,4383.173)
\texture{88555555 55000000 555555 55000000 555555 55000000 555555 55000000 
	555555 55000000 555555 55000000 555555 55000000 555555 55000000 
	555555 55000000 555555 55000000 555555 55000000 555555 55000000 
	555555 55000000 555555 55000000 555555 55000000 555555 55000000 }
\shade\path(3168,3713)(6243,5588)(6318,5513)
	(3243,3638)(3168,3713)
\path(3168,3713)(6243,5588)(6318,5513)
	(3243,3638)(3168,3713)
\path(4887,1913)(4212,113)
\dottedline{150}(4212,113)(2112,113)(12,2213)
	(12,4313)(4212,4313)(4212,113)
\put(4512,188){\makebox(0,0)[lb]{{\SetFigFont{7}{8.4}{\rmdefault}{\mddefault}{\updefault}$y_0$}}}
\put(4512,2063){\makebox(0,0)[lb]{{\SetFigFont{7}{8.4}{\rmdefault}{\mddefault}{\updefault}$x_2$}}}
\put(312,2138){\makebox(0,0)[lb]{{\SetFigFont{7}{8.4}{\rmdefault}{\mddefault}{\updefault}$Q = P \cap L$}}}
\put(6837,4538){\makebox(0,0)[lb]{{\SetFigFont{7}{8.4}{\rmdefault}{\mddefault}{\updefault}$\F_{x_1} \cap \F_{x_2}$}}}
\put(2412,4388){\makebox(0,0)[lb]{{\SetFigFont{7}{8.4}{\rmdefault}{\mddefault}{\updefault}$x_1$}}}
\end{picture}
}
  \caption{$\V(P) \subset (\F_{x_1} \cap \F_{x_2}) \cup L$}
  \label{fig:allInL}
\end{figure}
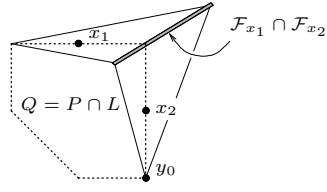
\item
\pagebreak
\underline{$\V(P) \not\subset (\F_{x_1} \cap \F_{x_2}) \cup L$ :} 
There are two subcases.
\begin{itemize}
\item[$\circ$]
\underline{$\V(P) \subset \F_{x_1} \cup \F_{x_2} \cup L$ :} 
Since $y_0 \in \V(\F_{x_2}) \cap L$, 
there has to exist a pair of vertices $y_1 \in \F_{x_1}\backslash\F_{x_2}$ and 
$y_2 \in \F_{x_2}\backslash\F_{x_1}$ with $y_1 + y_2 \not= 0$, such that 
one of the two sets $\{x_1,y_1,y_2\}$ or $\{x_2,y_1,y_2\}$ is linearly independent, say
the latter.

Choose $u \in \NR$ so that $\pro{u}{x_2} = \pro{u}{y_2} = 0$ and
$\pro{u}{y_1} = 1$. 
Set $A := \{r \in \F_{x_1} \cap \J \setminus \Lambda \;:\;
\pro{\eta_{x_2}}{r} = 0\}$.
Since $x_1 \in A$, $A \not= \emptyset$, so there exists $r \in A$ with
$\pro{u}{r}$ maximal (in particular $\ge 0$).

By Corollary~\ref{tool} we have $k_1 := \pro{\eta_{x_2}}{y_1} \geq 1$
and $k_2 := \pro{\eta_{x_1}}{y_2}$\\$\geq 1$.
Furthermore $q := p(r,y_2) = k_2 r + y_2 \in \J \setminus \Lambda$,
$\pro{\eta_{x_1}}{q} = 0$.
Since $\pro{\eta_{x_2}}{r} = 0$ and $y_2 \in \F_{x_2}$, we get $q \in
\F_{x_2}$.
Again by Corollary~\ref{tool}
$r' := p(q,y_1) = k_1 q + y_1 \in \J \setminus \Lambda$,
$\pro{\eta_{x_2}}{r'} = 0$. 
Since $\pro{\eta_{x_1}}{q} = 0$ and $y_1 \in \F_{x_1}$, we get $r' \in
\F_{x_1}$. Hence $r' \in A$. 
Since $r' = k_1 (k_2 r + y_2) + y_1$, we have 
$\pro{u}{r'} = k_1 k_2 \pro{u}{r} + k_1 > \pro{u}{r}$, a contradiction.
\begin{figure}[htbp]
  \centering
  \setlength{\unitlength}{0.00016667in}
\begingroup\makeatletter\ifx\SetFigFont\undefined%
\gdef\SetFigFont#1#2#3#4#5{%
  \reset@font\fontsize{#1}{#2pt}%
  \fontfamily{#3}\fontseries{#4}\fontshape{#5}%
  \selectfont}%
\fi\endgroup%
{\renewcommand{\dashlinestretch}{30}
\begin{picture}(7974,5889)(0,-10)
\put(3987,4812){\blacken\ellipse{212}{212}}
\put(3987,4812){\ellipse{212}{212}}
\put(3012,4212){\blacken\ellipse{212}{212}}
\put(3012,4212){\ellipse{212}{212}}
\put(6237,2712){\blacken\ellipse{212}{212}}
\put(6237,2712){\ellipse{212}{212}}
\put(6237,612){\blacken\ellipse{212}{212}}
\put(6237,612){\ellipse{212}{212}}
\put(837,4212){\blacken\ellipse{212}{212}}
\put(837,4212){\ellipse{212}{212}}
\path(12,3612)(4212,3612)(4212,12)
\path(3762,5862)(7962,5862)(7962,3237)
\path(4212,3612)(7962,5862)
\path(7962,2412)(7962,1512)
\dottedline{150}(3987,4812)(6237,4812)(6237,612)
\dottedline{150}(837,4212)(5187,4212)(6237,2712)
\put(5487,87){\makebox(0,0)[lb]{\smash{{{\SetFigFont{7}{8.4}{\rmdefault}{\mddefault}{\updefault}$y_2$}}}}}
\put(6462,2562){\makebox(0,0)[lb]{\smash{{{\SetFigFont{7}{8.4}{\rmdefault}{\mddefault}{\updefault}$q=p(r,y_2)$}}}}}
\put(1512,4512){\makebox(0,0)[rb]{\smash{{{\SetFigFont{7}{8.4}{\rmdefault}{\mddefault}{\updefault}$y_1$}}}}}
\put(3087,4812){\makebox(0,0)[lb]{\smash{{{\SetFigFont{7}{8.4}{\rmdefault}{\mddefault}{\updefault}$r$}}}}}
\put(2262,4362){\makebox(0,0)[lb]{\smash{{{\SetFigFont{7}{8.4}{\rmdefault}{\mddefault}{\updefault}$r'$}}}}}
\put(4512,1812){\makebox(0,0)[lb]{\smash{{{\SetFigFont{7}{8.4}{\rmdefault}{\mddefault}{\updefault}$\F_{x_2}$}}}}}
\put(4812,5187){\makebox(0,0)[lb]{\smash{{{\SetFigFont{7}{8.4}{\rmdefault}{\mddefault}{\updefault}$\F_{x_1}$}}}}}
\end{picture}
}
  \caption{Ping-pong}
  \label{fig:pingpong}
\end{figure}
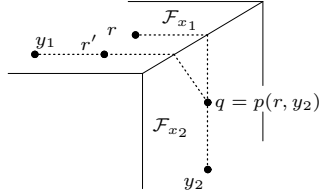
\item[$\circ$]
\underline{$\V(P) \not\subset \F_{x_1} \cup \F_{x_2} \cup L$ :} Let 
$y_1 \in \V(P)$ outside $\F_{x_1} \cup \F_{x_2} \cup L$. Now $x_3 := p(x_1,y_1) \in \J
\setminus \Lambda$ does not belong to $\F_{x_2}$ or $L$.
By Lemma~\ref{lemma:nill}, $z := z(x_2,x_3)=ax_2+bx_3 \in \randp \cap
\Lambda$ with $\pro{x_1}{z} = b \pro{x_1}{x_3} = 0$, a contradiction to Corollary~\ref{tool}.
\begin{figure}[htbp]
  \centering
  \setlength{\unitlength}{0.00016667in}
\begingroup\makeatletter\ifx\SetFigFont\undefined%
\gdef\SetFigFont#1#2#3#4#5{%
  \reset@font\fontsize{#1}{#2pt}%
  \fontfamily{#3}\fontseries{#4}\fontshape{#5}%
  \selectfont}%
\fi\endgroup%
{\renewcommand{\dashlinestretch}{30}
\begin{picture}(9087,5540)(0,-10)
\put(4275,1913){\blacken\ellipse{212}{212}}
\put(4275,1913){\ellipse{212}{212}}
\put(6075,1913){\blacken\ellipse{212}{212}}
\put(6075,1913){\ellipse{212}{212}}
\put(4275,113){\blacken\ellipse{212}{212}}
\put(4275,113){\ellipse{212}{212}}
\put(8325,3263){\blacken\ellipse{212}{212}}
\put(8325,3263){\ellipse{212}{212}}
\put(5400,4388){\blacken\ellipse{212}{212}}
\put(5400,4388){\ellipse{212}{212}}
\path(1875,3713)(6075,3713)(6075,113)
\path(4875,5513)(9075,5513)(9075,1613)
\path(6075,3713)(9075,5513)
\dottedline{150}(4275,1913)(6075,1913)(8325,3263)
\dottedline{150}(4275,113)(4275,3713)(5400,4388)
\put(4425,263){\makebox(0,0)[lb]{\smash{{{\SetFigFont{7}{8.4}{\rmdefault}{\mddefault}{\updefault}$y_1$}}}}}
\put(5175,4613){\makebox(0,0)[lb]{\smash{{{\SetFigFont{7}{8.4}{\rmdefault}{\mddefault}{\updefault}$x_1$}}}}}
\put(8025,3488){\makebox(0,0)[lb]{\smash{{{\SetFigFont{7}{8.4}{\rmdefault}{\mddefault}{\updefault}$x_2$}}}}}
\put(0,1763){\makebox(0,0)[lb]{\smash{{{\SetFigFont{7}{8.4}{\rmdefault}{\mddefault}{\updefault}$x_3=p(x_1,y_1)$}}}}}
\put(5475,2063){\makebox(0,0)[lb]{\smash{{{\SetFigFont{7}{8.4}{\rmdefault}{\mddefault}{\updefault}$z$}}}}}
\end{picture}
}
  \caption{Finish}
  \label{fig:ortho}
\end{figure}
\end{itemize}
\end{itemize}
\end{proof}

\end{document}